\documentclass{birkjour}
\usepackage[dvips]{epsfig}
\usepackage{graphics}
\usepackage{amsmath}
\usepackage[full]{textcomp}
\usepackage{amsfonts}
\usepackage{amsmath}
\usepackage{amsfonts}
\usepackage[psamsfonts]{amssymb}
\usepackage{amsbsy}
\usepackage{amsxtra}
\usepackage{amsthm}
\usepackage{mathabx}
\usepackage{cmll}
\usepackage{esint}
\usepackage[mathscr]{eucal}
\usepackage[cp1251]{inputenc}
\usepackage[OT1,T1,TS1,T2A]{fontenc}
\usepackage[russian,english]{babel}

\newtheorem*{nonumtheorem}{Theorem}

\newtheorem*{nonumdefinition}{Definition}

\newtheorem{theorem}{Theorem}
\newtheorem{lemma}{Lemma}[section]
\newtheorem{remark}{Remark}[section]

\def\theequation{\arabic{section}.\arabic{equation}}
  
\numberwithin{equation}{section}

\DeclareMathOperator{\re}{Re}
\DeclareMathOperator{\im}{Im}
\DeclareMathOperator{\ext}{Ext}
\DeclareMathOperator{\inte}{Int}
\DeclareMathOperator{\grad}{grad}

\begin{document}

%
%
%
%
%
\title{The function \(\boldsymbol{\cosh\big(\sqrt{a\,t^2+b}\big)}\) is exponentially convex.}
\author[V. Katsnelson]{Victor Katsnelson}
\address{%
Department of Mathematics\\
The Weizmann Institute\\
76100, Rehovot\\
Israel}
\email{victor.katsnelson@weizmann.ac.il; victorkatsnelson@gmail.com}

\begin{abstract} Given positive numbers \(a\) and \(b\), the function \(\sqrt{at^2+b}\)
is exponentially convex function of \(t\) on the whole real axis. Three proofs
of this result are presented.
\end{abstract}
\maketitle
\section{The exponential convexity result}
 \begin{nonumdefinition}
A function \(f\) on \(\mathbb{R}\), \(f:\,\mathbb{R}\to[0,\infty)\), is said to be \emph{exponentially convex}
if
\begin{enumerate}
\item[\textup{1}.]
 For every positive integer \(N\), for every choice of real numbers \(t_1,t_2,\,\ldots\,\), \(t_{N}\), and complex numbers
\(\xi_1\), \(\xi_2, \,\ldots\,, \xi_{N}\), the inequality holds
\begin{equation}
\label{pqf}
\sum\limits_{r,s=1}^{N}f(t_r+t_s)\xi_r\overline{\xi_s}\geq 0;
\end{equation}
\item[\textup{2}.]
The function \(f\) is continuous on \(\mathbb{R}\).
\end{enumerate}
\end{nonumdefinition}
The class of exponentially convex functions was introduced by S.N.Bernstein, \cite{B}, see \S 15 there.  Russian translation of the paper \cite{B} can be found in \cite[pp.370-425]{B1}.

From \eqref{pqf} it follows that
the inequality
\begin{math}
f(t_1+t_2)\leq\sqrt{f(2t_1)f(2t_2)}
\end{math}
holds for every \(t_1\in\mathbb{R},t_2\in\mathbb{R}\). Thus the alternative takes place: \\                        \textit{If \(f\) is an exponentially convex function then either \(f(t)\equiv 0\), or \(f(t)>0\) for every \(t\in\mathbb{R}\).}

\newpage
\noindent
\begin{center}
\textbf{Properties of the class of exponentially convex functions.}
\end{center}
\begin{enumerate}
\item[\textup{P\,1}.] If \(f(t\) if an exponentially convex function and \(c\geq0\) is a nonnegative constant, then the function \(cf(t)\) is exponentially convex.
\item[\textup{P\,2}.] If \(f_1(t)\) and \(f_2(t)\) are exponentially convex functions, then their sum
\(f_1(t)+f_2(t)\) is exponentially convex.
\item[\textup{P\,3}.] If \(f_1(t)\) and \(f_2(t)\) are exponentially convex functions, then their product
\(f_1(t)\cdot f_2(t)\) is exponentially convex.
\item[\textup{P\,4}.] Let \(\lbrace f_{n}(t)\rbrace_{1\leq n<\infty}\) be a sequence of exponentially
convex functions. We assume that for each \(t\in\mathbb{R}\) there exists the limit
\(f(t)=\lim_{n\to\infty}f_{n}(t)\), and that \(f(t)<\infty\ \forall t\in\mathbb{R}\).
Then the limiting function \(f(t)\) is exponentially convex.

\end{enumerate}

From the functional equation for the exponential function
it follows that for each real number \(\lambda\), for every choice of real numbers \(t_1,t_2,\,\ldots\,\), \(t_{N}\) and complex numbers
\(\xi_1\), \(\xi_2, \,\ldots\,, \xi_{N}\), the equality holds
\begin{equation}
\label{ece}
\sum\limits_{r,s=1}^{N}e^{\lambda(t_r+t_s)}\xi_r\overline{\xi_s}=
\bigg|\sum\limits_{p=1}^{N}e^{\lambda t_p}\xi_p\,\bigg|^{\,2}\geq 0.
\end{equation}
The relation \eqref{ece} can be formulated as
\begin{lemma}
\label{ECE}
For each real \(\lambda\), the function \(e^{\lambda t}\) of the variable \(t\) is exponentially convex.
\end{lemma}
For \(z\in\mathbb{C}\), the function \(\cosh z\), which is called \emph{the hyperbolic cosine of \(z\)}, is defined as
\begin{equation}
\label{dch}
\cosh z =\frac{1}{2}\big(e^z+e^{-z}\big).
\end{equation}
From Lemma \ref{ECE} and property P\,2 we obtain
\begin{lemma}
\label{echc}
For each real \(\mu\), the function \(\cosh(\mu\, t)\) of the variable \(t\) is exponentially convex.
\end{lemma}
The following result is well known.
\begin{nonumtheorem}[The representation theorem] {\ }\\[-4.0ex]
\begin{enumerate}
\item[\textup{1}.]
Let \(\sigma(d\lambda)\) be a nonnegative measure on the real axis,
and let the function \(f(t)\) be a two-sided Laplace transform of the measure
 \(d\sigma(\lambda)\):
\begin{equation}
\label{rep}
f(t)=\int\limits_{\lambda\in\mathbb{R}}e^{\lambda t}\,d\sigma(\lambda)
\end{equation}
for any \(t\in\mathbb{R}\). Then the function \(f\) is exponentially convex.
 \item[\textup{2}. ] Let \(f(t)\) be an exponentially convex function. Then this function \(f\) can be
  represented on \(\mathbb{R}\)  as a two-sided Laplace transform \eqref{rep} of a nonnegative measure \(d\sigma(\lambda)\).  \textup{(}In particular, the integral in the right
     hand side of \eqref{rep} is finite for any \(t\in\mathbb{R}\).\textup{)} The representing measure \(d\sigma(\lambda)\) is unique.
 \end{enumerate}
\end{nonumtheorem}

The assertion 1 of the representation theorem is an evident consequence of Lemma~\,\ref{ECE}, of the properties P\,1,\,P\,2, P\,4, and of the definition of the integration.

 The proof of the assertion 2 can be found in \cite{A},\,Theorem 5.5.4, and in
\cite{Wi},\,Theorem 21.

Of course, Lemma \ref{echc} is a special case of the representation theorem which corresponds to the representing measure
\(\sigma(d\lambda)=1/2\big(\delta(\lambda-\mu)+\delta(\lambda+\mu)\big)\,d\lambda\), where \(\delta(\lambda\mp\mu)\) are Dirak's \(\delta\)-functions supported at the points \(\pm\mu\).\\

\begin{lemma}
\label{ef}
The expression
\begin{equation}
 \label{mcf}
 \varphi(t,a,b)=\cosh\big(\sqrt{at^2+b}\,\big).
 \end{equation}
 is well defined for every complex numbers \(t,a,b\).
The function \(\varphi(t,a,b)\)
 is an entire function of complex variables \((t,a,b)\in\mathbb{C}^3\).  For each fixed \(a>0\) and \(b\),
the function \(\varphi(t,a,b)\), considered as a function of \(t\), is an entire function
of exponential type  \(\sqrt{a}\).
\end{lemma}
\noindent
\textsl{\textsf{Proof}}. The function \(\varphi(t,a,b)\)
is a superposition of the entire function \(\cosh\sqrt{\zeta}\) of variable \(\zeta\) and the quadratic polynomial \(\zeta(t,a,b)=at^2+b\). The assertion concerning the growth of this function is evident.
\hspace*\fill\(\Box\)

In the paper \cite{BMV} a conjecture was
formulated which now is commonly known as the BMV conjecture:\\[1.0ex]
\textbf{The BMV Conjecture.} Let \(A\) and \(B\) be Hermitian matrices of size
\(n\times{}n\).
Then the function
\begin{equation}
\label{TrF}
f_{A,B}(t)=
\textup{trace}\,\{\exp[tA+B]\}
\end{equation}
of the variable \(t\) is representable as a bilateral Laplace transform of a \textsf{non-negative} measure
\(d\sigma_{A,B}(\lambda)\) compactly supported on the real axis:
\begin{equation}
\label{LaR}
f_{A,B}(t)=\!\!\int\limits_{\lambda\in(-\infty,\infty)}\!\!\exp(t\lambda)\,d\sigma_{A,B}(\lambda), \ \ \forall \,t\in(-\infty,\infty).
\end{equation}

In general case, if the matrices \(A\) and \(B\) do not commute, the BMV conjecture remained an
open question for longer than 35 years. In 2011, Herbert
Stahl, \cite{St}, gave an affirmative answer to the BMV conjecture.\\[1.0ex]
\textbf{Theorem}\,(H.Stahl) \textit{Let \(A\) and \(B\) be \(n\times{}n\) Hermitian matrices.}
\textit{Then the function \(f_{A,B}(t)\) defined by \eqref{TrF} is representable
as the bilateral Laplace transform \eqref{LaR} of a non-negative measure \(d\sigma_{A,B}(\lambda)\)
supported on the closed interval \([\lambda_{\min},\lambda_{\max}]\).}

The proof of Herbert Stahl is based on ingenious considerations related to Riemann surfaces of algebraic functions.

In the case \(2\times 2\) matrices \(A\) and \(B\), the BMV conjecture is equivalent to the exponential convexity with respect to \(t\) for each
 \(a>0\) and \(b\geq 0\) of the function
\(\varphi(t,a,b)=\cosh\big(\sqrt{at^2+b}\,\big)\) which was introduced in \eqref{mcf} .

\begin{theorem}
 \label{MT}
 For each fixed \(a>0\) and \(b\geq0\), the function \(\varphi(t,a,b)\) defined by \eqref{mcf}
is an exponentially convex function of variable \(t\).
\end{theorem}

In what follow we present  three different proofs of Theorem \ref{MT}.  The first and the second proofs are based on the representation theorem. We prove that the function \(\widehat{d}(\lambda,a,b)\) which is defined by \eqref{rfv} below takes positive values
for \(\lambda\in(-\sqrt{a},\sqrt{a})\). In the first proof we calculate the function
 \(\widehat{d}(\lambda,a,b)\) explicitly expressing this function in terms of the modified
Bessel function \(I_{1}\). In the second proof, we prove the positivity of the function
\(\widehat{d}(\lambda,a,b)\) using the reasoning by Herbert Stahl in \cite{St}.
(We use a very simple special case of this reasoning.)
The third proof
is based on the Taylor expansion of the function \(\varphi(t,a,b)\), \eqref{mcf}, with respect to parameter \(b\). This proof does not use any integration in the complex plane. It based only on Lemma \ref{echc} and on the properties P1 -- P4
of the class of exponentially convex functions.
As a by-product of this proof we obtain that all coefficients of this Taylor expansion  are
exponentially convex functions. However we can not conclude directly  from this proof  that the
restriction of the representing measure on the \emph{open} interval \((-\sqrt{a},\sqrt{a})\)
is an absolutely continuous measure.

\begin{lemma}
\label{repl}
For each fixed \(a>0,\,b\geq0\), the function \(\varphi(t,a,b)\) defined by \eqref{mcf} is representable in the form
\begin{equation}
\label{rfv}
\varphi(t,a,b)=\cosh \sqrt{a}\,t+
\int\limits_{-\sqrt{a}}^{\sqrt{a}}\widehat{d}(\lambda,a,b)e^{\lambda t}\,d\lambda,
\ \ \ \forall t\in\mathbb{C},
\end{equation}
where the function \(\widehat{d}(\lambda,a,b)\) possesses the properties\\
\begin{enumerate}
\item[\textup{1.}]  {\ }
\vspace{-4.0ex}
\begin{equation}
\label{L2}
 \int\nolimits_{-\sqrt{a}}^{\sqrt{a}}|\widehat{d}(\lambda,a,b)|^2d\lambda<\infty;
 \end{equation}
 \item[\textup{2.}] The function  \(\widehat{d}(\lambda,a,b)\)
 is continuous with respect to  \(\lambda\) on the closed interval \([-\sqrt{a},\sqrt{a}]\), takes real values there, and is even.
\item[\textup{3.}] The values of the function \(\widehat{d}(\lambda,a,b)\) at the end points \(\pm\sqrt{a}\)
of the interval \([-\sqrt{a},\sqrt{a}]\) are:
     \begin{equation}
     \label{vep}
     \widehat{d}(\pm\sqrt{a},a,b)=\frac{b}{4\sqrt{a}}.
     \end{equation}
\end{enumerate}
\end{lemma}

\vspace{2.0ex}
\noindent
\textsl{Proof.} We introduce the function
\begin{equation}
\label{df}
d(t,a,b)=\cosh \big(\sqrt{at^2+b}\big)-\cosh \sqrt{a}t
\end{equation}
of variables \(t,a,b\). Considered as a function of \(t\) for fixed \emph{positive} \(a\) and \(b\),
\(d(t,a,b)\) is entire function of exponential type \(\sqrt{a}\).
On the imaginary axis \(d\) takes the form
\begin{equation}
\label{dia}
d(i\tau,a,b)=\cos\sqrt{a\tau^2-b}-\cos\sqrt{a}\tau, \quad \tau\in\mathbb{R}.
\end{equation}
From \eqref{dia} it follows that the function \(d\) is a bounded and decaying on the imaginary axis:
\(|d(i\tau,a,b)\leq{}1+\cosh b,\,\tau\in\mathbb{R}\), \(d(i\tau,a,b)=O(|\tau|^{-1}) \ \text{ as } \ \tau\to\pm\infty\).
By the Wiener-Paley theorem, the function \(d(i\tau,a,b)\) is representable in the form
\begin{equation}
\label{pwr}
d(i\tau,a,b)=\int\limits_{-\sqrt{a}}^{\sqrt{a}}\widehat{d}(\lambda,a,b)e^{i\lambda \tau}\,d\lambda, \ \
\tau\in\mathbb{R},
\end{equation}
where the function \(\widehat{d}(\lambda,a,b)\) satisfies the condition \eqref{L2}.
\emph{The equality \eqref{pwr} serves as a definition of the function \(\widehat{d}(\lambda,a,b)\).}
So, this function is defined only for \(a>0,\,b>0,\,-\sqrt{a}\leq\lambda\leq\sqrt{a}\).

Since the function \(d(i\tau,a,b)\) is even with respect to \(\tau\) and real valued, its inverse Fourier
transform \(\widehat{d}(\lambda,a,b)\) is even with respect to \(\lambda\) and real valued.

From \eqref{dia} we obtain that
\begin{equation}
\label{ase}
d(i\tau,a,b)-\frac{b}{2}\,\frac{\sin\sqrt{a}\,\tau}{\sqrt{a}\,\tau}=O(\tau^{-2}) \text{ as } \tau\to\pm\infty.
\end{equation}
Hence the function in the left hand side of \eqref{ase} is a Fourier transform of some function
 \(r(\lambda)\) which is square summable and \emph{continuous} at every \(\lambda\in\mathbb{R}\).
We remark that
\begin{equation*}
\frac{b}{2}\,\frac{\sin\sqrt{a}\,\tau}{\sqrt{a}\,\tau}=
\int\limits_{-\sqrt{a}}^{\sqrt{a}}\widehat{c}\,e^{i\lambda\tau}\,d\lambda, \ \ \tau\in\mathbb{R},
\end{equation*}
where \(\widehat{c}=\frac{b}{4\sqrt{a}}\) is a constant function. Hence
\begin{equation*}
r(\lambda)=
\begin{cases}
\widehat{d}(\lambda,a,b)-\widehat{c},& \ \ \textup{for} \ |\lambda|<\sqrt{a},\\
0\phantom{( \ \ \ \ \     \lambda)-\widehat{c}},& \ \ \textup{for} \ |\lambda|>\sqrt{a}.
\end{cases}
\end{equation*}
 Since \(r(\lambda)=0\) for \(|\lambda|>\sqrt{a}\),
also \(r(\pm\sqrt{a})=0\). Thus, \eqref{vep} holds.          \hspace\fill \(\Box\)

\vspace{2.0ex}

\section{Representation the function \(\boldsymbol{\widehat{d}(\lambda,a,b)}\)\\ by a contour integral.\label{RCI}}
Let \(S\) be a segment of the imaginary axis:
\begin{equation}
\label{slit}
S=\left\lbrace\zeta=%
\xi+i\eta:\,\xi=0,-\sqrt{\tfrac{b}{a}}\leq\eta\leq\sqrt{\tfrac{b}{a}}\,\right\rbrace.
\end{equation}
The function \(\sqrt{a\zeta^{2}+b}\) is a single value function of \(\zeta\) in the complex plane
slitted along the vertical segment \(S\).
We choose the branch of this function which takes positive values for large real \(\zeta\).
\begin{lemma}
\label{lir}
The function \(\widehat{d}(\lambda,a,b)\), which was defined by \eqref{pwr}, admits the integral representation
\begin{equation}
\label{cir}
\widehat{d}(\lambda,a,b)=-\frac{1}{4\pi i}\ointctrclockwise\limits_{\Gamma}e^{-\sqrt{a\zeta^2+b}}%
e^{-\lambda\zeta}\,d\zeta, \ \
-\sqrt{a}<\lambda<\sqrt{a},
\end{equation}
where \(\Gamma\) is an arbitrary counterclockwise oriented closed Jordan curve which contains the slit \(S\)  inside.
\end{lemma}
\noindent
\textsl{Proof.} According the inversion formula for the Fourier transform,
\begin{equation}
\label{iff}
\widehat{d}(\lambda,a,b)=
\frac{1}{2\pi}\int\limits_{-\infty}^{\infty}\varphi(i\eta,a,b)e^{-i\lambda\eta}\,d\eta.
\end{equation}
We interpret the integral in the right hand side of \eqref{iff} as the integral along the the vertical
straight line \(\lbrace\zeta:\, \re\zeta =0\rbrace\):
\begin{equation}
\label{iffv}
\widehat{d}(\lambda,a,b)=\frac{1}{2\pi i}\int\limits_{\re\zeta=0}\varphi(\zeta,a,b)e^{-\lambda\zeta}\,d\zeta
=\lim_{R\to+\infty}\frac{1}{2\pi i}\int\limits_{-iR}^{+iR}\varphi(\zeta,a,b)e^{-\lambda\zeta}\,d\zeta.
\end{equation}
Since the function \(\varphi(\zeta,a,b)\) is bounded in each vertical strip
\(\lbrace\zeta:\,\alpha\leq\re\zeta\leq\beta\rbrace\) and tends to zero as \(\im\zeta\to\pm\infty\) within this strip, the value of the integral in \eqref{iffv} does not change if we integrate
along any vertical line \(\lbrace\zeta:\,\re\zeta=\gamma\rbrace\), where
\(\gamma\) is an arbitrary real number:
\begin{equation}
\label{iffvg}
\widehat{d}(\lambda,a,b)=\frac{1}{2\pi i}\int\limits_{\re\zeta=
\gamma}\varphi(\zeta,a,b)e^{-\lambda\zeta}\,d\zeta,\ \ -\sqrt{a}\leq\lambda\leq\sqrt{a}\,.
\end{equation}
Choosing \(\gamma<0\) , we split the integral in \eqref{iffvg} into the sum
\begin{equation}
\label{iffsp}
\widehat{d}(\lambda,a,b)=\frac{1}{2\pi i}\int\limits_{\re\zeta=\gamma}\varphi_{+}(\zeta,a,b)e^{-\lambda\zeta}\,d\zeta+
\frac{1}{2\pi i}\int\limits_{\re\zeta=\gamma}\varphi_{-}(\zeta,a,b)e^{-\lambda\zeta}\,d\zeta,
\end{equation}
where
\begin{equation}
\varphi_{+}(\zeta,a,b)=\frac{1}{2}\Big(e^{\sqrt{a\zeta^{2}+b}}-e^{\sqrt{a}\zeta}\Big),\quad
\varphi_{-}(\zeta,a,b)=\frac{1}{2}\Big(e^{-\sqrt{a\zeta^{2}+b}}-e^{-\sqrt{a}\zeta}\Big).
\label{ppm}
\end{equation}

The function \(\varphi_{+}(\zeta,a,b)\) is holomorphic in the halfplane \(\lbrace\zeta:\,\re\zeta\leq\gamma\rbrace\) and admits the estimate
\begin{equation*}
|\varphi_{+}(\zeta,a,b)|\leq  c(\gamma)(1+|\zeta|)^{-1}e^{\sqrt{a}\re\zeta}, \quad \forall \zeta: \re\zeta\leq\gamma
\end{equation*}
there, where \(c(\gamma)<\infty \) is a constant. Therefore
\begin{equation}
\label{ipi+}
\frac{1}{2\pi i}\int\limits_{\re\zeta=\gamma}\varphi_{+}(\zeta,a,b)e^{-\lambda\zeta}\,d\zeta=0,
\ \ \lambda<\sqrt{a} .
\end{equation}
The function \(\varphi_{-}(\zeta,a,b)\) is holomorphic in the slitted half plane
\(\lbrace\zeta:\,\re\zeta\geq\gamma,\)
\(\zeta\not\in S\rbrace\) and admits the
estimate
\begin{equation*}
|\varphi_{-}(\zeta,a,b)|\leq  c(\gamma)(1+|\zeta|)^{-1}e^{-\sqrt{a}\re\zeta}, \quad \forall \zeta: \re\zeta\geq\gamma,\,\zeta\not\in S
\end{equation*}
there. Therefore
\begin{equation}
\label{ipi-}
\frac{1}{2\pi i}\int\limits_{\re\zeta=\gamma}\varphi_{-}(\zeta,a,b)e^{-\lambda\zeta}\,d\zeta=
\frac{1}{2\pi i}\ointclockwise\limits_{\Gamma}\varphi_{-}(\zeta,a,b)e^{-\lambda\zeta}e^{-\lambda\zeta}\,d\zeta,
\ \ \lambda>-\sqrt{a} ,
\end{equation}
where \(\Gamma\) is an arbitrary closed Jourdan curve which is oriented \emph{clockwise} and contains
the slit \(S\) in its interior.
Since the function \(e^{-\sqrt{a}\zeta}\) is entire,
\(\oint\limits_{\Gamma}e^{-\sqrt{a}\zeta}\,d\zeta=0\). So
\begin{equation}
\label{red}
\frac{1}{2\pi i}\ointclockwise\limits_{\Gamma}\varphi_{-}(\zeta,a,b)e^{-\lambda\zeta}\,d\zeta=
-\frac{1}{4\pi i}\ointctrclockwise\limits_{\Gamma}e^{-\sqrt{a\zeta^2+b}}e^{-\lambda\zeta}\,d\zeta,
\end{equation}
where the integral in the right hand side is taken over the curve \(\Gamma\) which is oriented
\emph{counterclockwise}.
Comparing \eqref{iffsp}, \eqref{ipi+}, \eqref{ipi-}, and \eqref{red}, we obtain \eqref{cir}.
\hspace\fill \(\Box\)

\section{Explicite calculation of the function \(\boldsymbol{\widehat{d}(\lambda,a,b)}\).}
\begin{lemma}
\label{lsir}
The function \(\widehat{d}(\lambda,a,b)\) which was defined by \eqref{pwr} admits the integral representation
\begin{equation}
\label{ird}
\widehat{d}(\lambda,a,b)=\frac{1}{\pi}\sqrt{\frac{b}{a}}%
\int\limits_{0}^{1}\sinh\sqrt{b(1-\tau^2)}%
\cdot\cos\Big(\lambda\sqrt{\tfrac{b}{a}}\tau\Big)\,d\tau, \ \ -\sqrt{a}\leq\lambda\leq\sqrt{a}\,.
\end{equation}
\end{lemma}
\textsl{Proof.} We derive Lemma \ref{lsir} from Lemma \ref{lir} showing that
\begin{multline}
\label{dll}
-\frac{1}{4\pi i}\ointctrclockwise\limits_{\Gamma}e^{-\sqrt{a\zeta^2+b}}e^{-\lambda\zeta}\,d\zeta=
\frac{1}{\pi}\sqrt{\frac{b}{a}}\int\limits_{0}^{1}\sinh\sqrt{b(1-\tau^2)}%
\cdot\cos\lambda\sqrt{\tfrac{b}{a}}\tau\,d\tau.
\end{multline}
The function \(e^{-\sqrt{a\zeta^2+b}}e^{-\lambda\zeta}\) is holomorphic in the domain
 \(\mathbb{C}\setminus S\) and continuous up to boundary \(S=\partial(\mathbb{C}\setminus S)\) of this domain.
 Therefore the integral of this function over \(\Gamma\) does not change if we shrink the original contour \(\Gamma\) to the boundary \(S\):
 \begin{equation}
 \label{scr}
 \frac{1}{4\pi i}\ointctrclockwise\limits_{\Gamma}e^{-\sqrt{a\zeta^2+b}}e^{-\lambda\zeta}\,d\zeta=
 \frac{1}{4\pi i}\ointctrclockwise\limits_{S}e^{-\sqrt{a\zeta^2+b}}e^{-\lambda\zeta}\,d\zeta
 \end{equation}
 To one
 "geometric" point \(i\eta\in S\) there corresponds two topologically different "boundary"
 points \(+0+i\eta\) and \(-0+i\eta\) lying on the right edge \(S^{+}\) and the left edge
  \(S^{-}\) of the slit \(S\) respectively. The chosen branch of the function
  \(\sqrt{a\zeta^2+b}\) takes the following values on the boundary of the domain \(\mathbb{C}\setminus S\):
 \begin{equation}
 \label{onb}
 \sqrt{a(+0+i\eta)^2+b}=-\sqrt{a(-0+i\eta)^2+b}=\sqrt{b-a\eta^2}, \ \ i\eta\in S.
 \end{equation}
 If the point \(\zeta=\pm0+i\eta\) runs over \(S=\partial\, (\mathbb{C}\setminus S)\)
 counterclockwise, the \(\eta\) increases from \(-\sqrt{\tfrac{b}{a}}\) to \(\sqrt{\tfrac{b}{a}}\) if \(\zeta=\in S^{+}\) and  \(\eta\) decreases from \(\sqrt{\tfrac{b}{a}}\) to \(-\sqrt{\tfrac{b}{a}}\) if \(\zeta=\in S^{-}\).
 Therefore
 \begin{align*}
\sideset{}{\hspace*{-2.5ex}{\text{\raisebox{0.5ex}{\(\scriptstyle\curvearrowleft\)}}}}%
\int\limits_{S^{+}}e^{-\sqrt{a\zeta^2+b}}\,e^{-\lambda\zeta}\,d\zeta=&
+i\int\limits_{-\sqrt{b/a}}^{\sqrt{b/a}}%
e^{-\sqrt{b-a\eta^2}}e^{-i\lambda\eta}\,d\eta,\\
\sideset{}{\hspace*{-2.5ex}{\text{\raisebox{0.5ex}{\(\scriptstyle\curvearrowleft\)}}}}%
\int\limits_{S^{-}}e^{-\sqrt{a\zeta^2+b}}\,e^{-\lambda\zeta}\,d\zeta=&
-i\int\limits_{-\sqrt{b/a}}^{\sqrt{b/a}}%
e^{+\sqrt{b-a\eta^2}}e^{-i\lambda\eta}\,d\eta
\end{align*}
Thus
\begin{multline}
\label{ios}
\frac{1}{4\pi i}\ointctrclockwise\limits_{S}e^{-\sqrt{a\zeta^2+b}}e^{-\lambda\zeta}\,d\zeta=
\frac{1}{4\pi}\int\limits_{-\sqrt{b/a}}^{\sqrt{b/a}}%
(e^{-\sqrt{b-a\eta^2}}-e^{\sqrt{b-a\eta^2}})e^{-i\lambda\eta}\,d\eta=\\
=-\frac{1}{2\pi}\int\limits_{-\sqrt{b/a}}^{\sqrt{b/a}}%
\sinh\sqrt{b-a\eta^2}\,e^{-i\lambda\eta}\,d\eta=-\frac{1}{\pi}
\int\limits_{0}^{\sqrt{b/a}}\sinh\sqrt{b-a\eta^2}\,\cos\lambda\eta\,d\eta=\\
=-\frac{1}{\pi}\sqrt{\tfrac{b}{a}}\int\limits_{0}^{1}\sinh\sqrt{b(1-\eta^2)}
\cdot\cos\Big(\sqrt{\tfrac{b}{a}}\lambda\eta\Big)\,d\eta.
\end{multline}
Comparing \eqref{cir} with \eqref{ios}, we obtain \eqref{ird}.
\begin{lemma}
\label{lexex}
Let \(a>0\) and \(b>0\) be fixed positive numbers. Then
\begin{enumerate}
\item[\textup{1}.] The function \(\widehat{d}(\lambda,a,b)\)  which was defined by \eqref{pwr} can be expressed
explicitly in terms of the modified Bessel function \(I_1\):
\begin{equation}
\label{exex}
\widehat{d}(\lambda,a,b)=\frac{\sqrt{b}}{2\sqrt{a-\lambda^2}}
I_1\Big({\textstyle\sqrt{\frac{(a-\lambda^2)b}{a}}}\,\Big), \ \ \ \ -\sqrt{a}\leq\lambda\leq\sqrt{a}.
\end{equation}
\item[\textup{2}.]The function \(\widehat{d}(\lambda,a,b)\) is representable as the sum of the series
\begin{multline}
\label{dtb}
\widehat{d}(\lambda,a,b)=\frac{b}{4\sqrt{a}}\sum\limits_{k=0}^{\infty}\frac{1}{k!(k+1)!}
\bigg(\frac{(a-\lambda^2)b}{4a}\bigg)^k, \\ a>0, \ -\sqrt{a}\leq\lambda\leq\sqrt{a}, \ b\geq0.
\end{multline}
\end{enumerate}
\end{lemma}
\noindent

\begin{remark}
\label{dwrb}
The expression in the right hand sides of \eqref{ird} is an entire function of three variables
\((\lambda,\sqrt{a^{-1}},b)\in\mathbb{C}^3\). However the equalities \eqref{ird}, \eqref{exex},
\eqref{dtb} hold only for \(a>0,\,b>0,\, -\sqrt{a}\leq\lambda\leq\sqrt{a}\).
\textup{(}We recall that the function \(\widehat{d}(\lambda,a,b)\) was \emph{defined} by \eqref{pwr} only for \(a>0,\,b>0,\, -\sqrt{a}\leq\lambda\leq\sqrt{a}\).\textup{)}
\end{remark}

\noindent
\textsl{\textsf{Proof of Lemma \ref{lexex}.}}
We start from the formula \eqref{ird}. Using the Taylor expansion of the hyperbolic \(\sinh\)
function, we obtain
\begin{equation}
\label{ehs}
\widehat{d}(\lambda,a,b)=\frac{1}{\pi}\sqrt{\frac{b}{a}}\sum\limits_{r=0}^{\infty}\frac{1}{(2r+1)!}
b^{\,r+\frac{1}{2}}\int\limits_{0}^{1}(1-\tau^2)^{r+\frac{1}{2}}%
\cos\Big(\lambda\sqrt{\tfrac{b}{a}}\tau\big)\,d\tau
\end{equation}
The integral in the right hand side of \eqref{ehs} can be expressed in terms of the Bessel function
\(J_{r+1}\), see \cite[\textbf{9.1.20}]{AS}:
\begin{multline*}
\int\limits_{0}^{1}(1-\tau^2)^{r+\frac{1}{2}}%
\cos\Big(\lambda\sqrt{\tfrac{b}{a}}\tau\Big)\,d\tau=\\%
=\pi^{1/2}\,2^r\,\Gamma(r+3/2)a^{\frac{r+1}{2}}b^{-\frac{r+1}{2})}\,\lambda^{-(r+1)}\,
J_{r+1}\Big(\lambda\sqrt{\tfrac{b}{a}}\Big).
\end{multline*}
Substituting the last equality into \eqref{ehs}, we obtain the equality
\begin{equation*}
\widehat{d}(\lambda,a,b)=\pi^{-\frac{1}{2}}\sum\limits_{r=0}^{\infty}
\frac{\Gamma(r+3/2)}{(2r+1)!}2^ra^{\frac{r}{2}}b^{\frac{r+1}{2}}\lambda^{-(r+1)}
J_{r+1}\Big(\lambda\sqrt{\tfrac{b}{a}}\Big).
\end{equation*}
Taking into account the duplication formula for the Gamma-function, \cite[\textbf{6.1.18}]{AS}:
\begin{equation*}
\frac{\Gamma(r+\frac{3}{2})}{\Gamma(2r+2)}=\pi^{\frac{1}{2}}2^{-(2r+1)}\frac{1}{\Gamma(r+1)},
\end{equation*}
we transform the last equality to the form
\begin{equation}
\label{pmt}
\widehat{d}(\lambda,a,b)=\sum\limits_{r=0}^{\infty}\frac{1}{r!}2^{-(r+1)}
a^{\frac{r}{2}}b^{\frac{r+1}{2}}\lambda^{-(r+1)}
J_{r+1}\Big(\lambda\sqrt{\tfrac{b}{a}}\Big).
\end{equation}
Now we would like to reduce the equality \eqref{pmt} to the form which occurs in the so called
\emph{Multiplication Theorem}\footnote{
Proof of the Multiplication Theorem can be found in \cite[Chapter V, sec.5.22]{W}, see formula (15)
on page 142 of the English edition or on the page 156 of the Russian translation. See
also \cite[Chapter IV, sec.21]{Sc}.
}, see \cite[\textbf{9.1.74}]{AS}:
\begin{equation}
\label{omt}
\widehat{d}(\lambda)=\tfrac{\sqrt{b}}{2\lambda}\sum\limits_{r=0}^{\infty}\frac{(-1)^r}{r!}
\big(-\tfrac{a}{\lambda^2}\big)^r\cdot
\Big(\tfrac{\lambda}{2}\sqrt{\tfrac{b}{a}}\,\Big)^r\,J_{r+1}\Big(\lambda\sqrt{\tfrac{b}{a}}\Big).
\end{equation}
Let us introduce \(\mu: \mu^2-1=-\dfrac{a}{\lambda^2}\), i.e.
\begin{equation}
\label{lm}
\mu=i\frac{\sqrt{a-\lambda^2}}{\lambda}.
\end{equation}
Then the equality \eqref{omt} takes the form
\begin{equation}
\label{req}
\widehat{d}(\lambda,a,b)=\frac{\sqrt{b}}{2i\sqrt{a-\lambda^2}}\,\cdot\,\mu\sum\limits_{r=0}^{\infty}\frac{1}{r!}
(\mu^2-1)^r\Big(\tfrac{\lambda}{2}\sqrt{\tfrac{b}{a}}\,\Big)^r\,J_{r+1}\Big(\lambda\sqrt{\tfrac{b}{a}}\Big).
\end{equation}
According to the Multiplication Theorem,
\begin{equation}
\label{amt}
\mu\sum\limits_{r=0}^{\infty}\frac{1}{r!}
(\mu^2-1)^r\Big(\tfrac{\lambda}{2}\sqrt{\tfrac{b}{a}}\,\Big)^r\,
J_{r+1}\Big(\lambda\sqrt{\tfrac{b}{a}}\Big)=J_1\Big(i\sqrt{a-\lambda^2}\sqrt{\tfrac{b}{a}}\Big),
\ \ \forall \lambda,a,b.
\end{equation}
Taking into account that \(J_1(iz)=iI_1(z)\), we reduce the equality \eqref{req} to the form \eqref{exex}.

Using the Taylor expansion of the modified Bessel function \(I_1\), \cite[\textbf{9.6.10}]{AS},
we represent the function \(\widehat{d}(\lambda,a,b)\) as the sum of the series \eqref{dtb}.
\hspace*\fill\(\Box\)\\

\noindent
\section{The first proof of Theorem \ref{MT}.} From the equality \eqref{dtb} is evident that
\begin{equation}
\label{pod}
\widehat{d}(\lambda,a,b)>0 \text{ for } \lambda\in[-\sqrt{a},\sqrt{a}]
\end{equation}
Theorem \ref{MT} follows from \eqref{pod} and \eqref{rfv}. \hspace*\fill \(\Box\)\\

\noindent
\begin{theorem}
\label{ptc}
For each \(a>0\), the function \(\varphi(t,a,b)\) which was introduced in \eqref{mcf} admits the Taylor expansion with respect to \(b\):
\begin{subequations}
\label{tew}
\begin{equation}
\label{tewb}
\varphi(t,a,b)=\sum\limits_{k=0}^{\infty}\frac{1}{k!}\varphi_{k}(t,a)b^k, \ \ \forall\,t\in\mathbb{R},\
0\leq b<\infty.
\end{equation}
For each \(k\geq 0\), the function \(\varphi_{k}(t,a)\), which is the \(k\)-th coefficient of the Taylor expansion \eqref{tewb}, is exponentially convex:
\begin{align}
\varphi_0(t,a)=&\cosh\sqrt{a}t,\label{k=0}\\
\varphi_k(t,a)=&\frac{1}{(k-1)!4^ka^{k-\frac{1}{2}}}\int\limits_{-\sqrt{a}}^{\sqrt{a}}
(a-\lambda^2)^{k-1}\,e^{\lambda t}\,d\lambda,\ \ \ k=1,2,3,\,\ldots\,.
\label{k>0}
\end{align}
\end{subequations}
\end{theorem}
\noindent
\textsl{\textsf{Proof}.} The expansion \eqref{dtb} can be presented as a Taylor expansion with respect to
\(b\):
\begin{subequations}
\label{ted}
\begin{gather}
\widehat{d}(\lambda,a,b)=\sum\limits_{k=1}^{\infty}\frac{1}{k!}\widehat{d}_k(\lambda,a)\,b^k,
\label{ted1}\\
\intertext{where}
\widehat{d}_k(\lambda,a)=\frac{(a-\lambda^2)^{k-1}}{(k-1)!4^ka^{k-\frac{1}{2}}}, \ \ \ k=1,2,3,\,\ldots\,.
\label{ted2}
\end{gather}
\end{subequations}
Substituting the expansion \eqref{ted} into the integrand in \eqref{rfv}, we obtain the expansion
\eqref{tew}. It is evident that \(\widehat{d}_k(\lambda,a)>0\) for \(-\sqrt{a}<\lambda<\sqrt{a}\).
The exponential convexity of the function \(\varphi_k(t,a)\) follows from the
representation \eqref{k>0}.\hspace*\fill\(\Box\)\\

\begin{remark}
\label{rade}
The function \(\varphi_k(t,a)\), \eqref{k>0}, can be expressed in terms of the modified Bessel function
\(I_{k-\frac{1}{2}}\):
\begin{equation}
\label{ade}
\varphi_k(t,a)=\pi^{\frac{1}{2}}2^{-(k+\frac{1}{2})}a^{-\frac{k}{2}+\frac{1}{4}}
t^{-(k-\frac{1}{2})}\,I_{k-\frac{1}{2}}(\sqrt{a}t)\,,\\
\ \ a>0,\ t\in\mathbb{R},
  \ \ k=1,2,3,\,\ldots\,.
\end{equation}
See \textup{\cite[\textbf{9.6.18}]{AS}}.
\end{remark}

\begin{remark}
\label{prior1}
The formula \eqref{dtb} appeared in \textup{\cite[subsection 7.3]{St}}, see formulas \textup{(7.22)} and \textup{(7.23)} there. In \textup{\cite{St}}, the derivation of the expansion \eqref{dtb} was done by a direct calculation,
without any reference to the multiplication theorem for Bessel function. It should be mention that the
series in the right hand side of \eqref{dtb} appeared in \textup{\cite[section 2]{MK}} as a perturbation
series related to the BMV conjection for \(2\times2\) matrices.
\end{remark}

\section{The second proof of Theorem \ref{MT}.}
\noindent
 The starting point of the first as well as of the second is the representation of the value \(\widehat{d}(\lambda,a,b)\) by the contour integral \eqref{cir}. See Lemma \ref{lir}.

  In the first proof, we shrank the contour of integration over the slit \(S\),
  so the contour of integration was the same for every \(\lambda\in[-\sqrt{a},\sqrt{a}]\).

  In contrast to this, in the second proof we choose the contour \(\Gamma\) in such a way that the exponent \(-\sqrt{a\zeta^2+b}-\lambda\zeta\) of the integrand \(e^{-\sqrt{a\zeta^2+b}-\lambda\zeta}\) in \eqref{cir} takes real values on \(\Gamma\).
 (So the contour \(\Gamma\) depends on \(\lambda\)!). We denote this contour by \(\Gamma_{\lambda}\)

 The function \(\widehat{d}(\lambda,a,b)\) is even with
respect to \(\lambda\). Therefore
to prove the exponential convexity of the function \(\varphi(t,a,b)\), it is enough
to prove that the value \(\widehat{d}(\lambda,a,b)\) is positive for each
\begin{equation}
\label{negl}
\lambda\in(-\sqrt{a},0).
\end{equation}
 We choose an arbitrary \(\lambda\) satisfying the condition \eqref{negl} and fix this choice in the course of the proof.

Let us introduce the functions
\begin{subequations}
\label{Harf}
\begin{align}
\label{Harfr}
u(\zeta)=\re(\sqrt{a\zeta^2+b}+\lambda\zeta),\ \ z\in\mathbb{C}\setminus S, \\
v(\zeta)=\im(\sqrt{a\zeta^2+b}+\lambda\zeta),\ \ z\in\mathbb{C}\setminus S,
\label{Harfi}
\end{align}
\end{subequations}
where \(S\) is the vertical slit \eqref{slit} and the branch of the function \(\sqrt{a\zeta^2+b}\)
 in \(\mathbb{C}\setminus S\) is chosen which takes positive values for large real \(\zeta\).
 \begin{lemma}
\label{laze}
Let us assume that \(a>0,b>0\)  and \(\lambda\) satisfies the condition \eqref{negl}.
Then there exist \(\varepsilon>0\) \(R<\infty\), \(\varepsilon=\varepsilon(a,b,\lambda),\,R=R(a,b,\lambda)\), such that
\begin{subequations}
\label{poda}
\begin{align}
\label{podal}
v(\zeta)/\im\zeta>0,& \ \ \ \   \forall \zeta\in\mathbb{C}\phantom{\setminus S}:\,|\zeta|>R, \ \im\zeta\neq0,\\
\label{podas}
v(\zeta)/\im\zeta<0,& \ \ \ \   \forall \zeta\in\mathbb{C}\setminus S:\,|\zeta|<\varepsilon, \ \,\,\im\zeta\neq0.
\end{align}
\end{subequations}
\end{lemma}
\noindent
\textsf{\textsl{Proof.}} From the identity
\begin{equation*}
\sqrt{a\zeta^2+b}-\sqrt{a}\zeta=\frac{b}{\sqrt{a\zeta^2+b}+\sqrt{a}\zeta}
\end{equation*}
we derive that
\begin{equation*}
\im\sqrt{a\zeta^2+b}-\sqrt{a}\im\zeta=-\big(\im\sqrt{a\zeta^2+b}+\sqrt{a}\im\zeta\big)\,\rho(\zeta)
\end{equation*}
where \(\rho(\zeta)=b\big|\sqrt{a\zeta^2+b}+\sqrt{a}\zeta\big|^{-2}\). Thus
\begin{equation*}
\im\sqrt{a\zeta^2+b}=\frac{1-\rho(\zeta)}{1+\rho(\zeta)}\cdot\sqrt{a}\, \im\zeta
\end{equation*}
and
\begin{equation}
\label{exprv}
v(\zeta)=\bigg(\sqrt{a}\,\frac{1-\rho(\zeta)}{1+\rho(\zeta)}+\lambda\bigg)\cdot\im\zeta.
\end{equation}
It is clear that\footnote{\,Here the choice of the branch of the function \(\sqrt{a\zeta^2+b}\) is important.
}
\(\rho(\zeta)\to0\) as \(|\zeta|\to\infty\),  \(\rho(\zeta)\to1\) as \(|\zeta|\to0,\zeta\not\in S\).
Since \(\sqrt{a}+\lambda>0\), the inequality \eqref{podal} holds if if \(|\zeta|\) is large enough.
 Since \(\lambda<0\), the inequality \eqref{podas} holds if \(|\zeta|\) is small enough.
\hspace*\fill\(\Box\)\\[0.0ex]

Let \(N_{\lambda}\) be the set
 \begin{equation}
 \label{zers}
 N_{\lambda}=\lbrace\zeta\in\mathbb{C}\setminus S: v(\zeta)=0\rbrace,
 \end{equation}

 \begin{lemma} {\ }\\[-3.0ex]
 \label{szers}
 \begin{enumerate}
 \item[\textup{1}.]
 The set \(N_{\lambda}\) is the union of the real axis and an ellipse \(\Gamma_{\lambda}\)\textup{:}
 \begin{equation}
 \label{szs}
 N_{\lambda}=(\mathbb{R}\setminus0)+\Gamma_{\lambda}.
 \end{equation}
 where the ellipse \(\Gamma_{\lambda}\) is described by the equation:
 \begin{equation}
 \label{eqel}
 \frac{\xi^2}{A^2}+\frac{\eta^2}{B^2}=1, \ \ (\zeta=\xi+i\eta),
 \end{equation}
 with
 \begin{equation}
 \label{halfa}
 A=\sqrt{\frac{b}{a}}\,\cdot\,\frac{|\lambda|}{\sqrt{a}}\bigg(1-\frac{\lambda^2}{a}\bigg)^{-\frac{1}{2}},\quad
 B=\sqrt{\frac{b}{a}}\,\cdot\,\bigg(1-\frac{\lambda^2}{a}\bigg)^{-\frac{1}{2}}
 \end{equation}
 \item[\textup{2}.] The slit \(S\) is contained in the interior of the ellipse \(\Gamma_\lambda\).
 \end{enumerate}
 \end{lemma}
\noindent
\textsf{\textsl{Proof.}}\\
\textsf{1.} Let \(\zeta=\xi+i\eta,\,\,\sqrt{a\zeta^2+b}=p+iq\), where \(\xi,\eta,p,q\)
are real numbers. The equality
\begin{equation*}
\pm\sqrt{a\zeta^2+b}=p+iq
\end{equation*}
is equivalent to the system of equalities
\begin{equation}
\label{syst}
\begin{cases}
a(\xi^2-\eta^2)+b&=p^2-q^2,\\
\ \ \ \ a\xi\eta&=\ \ pq.
\end{cases}
\end{equation}
Here \(p=p(\xi,\eta),\,q=q(\xi,\eta).\) Clearly \(v(\xi,\eta)=q(\xi,\eta)+\lambda\eta\).

Let \(\zeta\in N_{\lambda}\). This means that \(v(\xi,\eta)=0\), i.e.
\begin{subequations}
\label{pq}
 \begin{equation}
 \label{q}
 q=-\lambda\eta
 \end{equation}
 Substituting this equality into the second equality of the system \eqref{syst},
we obtain the equality \(a\xi\eta=-\lambda p\eta.\) Assuming that \(\eta\neq0\), that is \(\zeta\not\in\mathbb{R}\),
we can cancel by \(\eta\) and obtain
\begin{equation}
\label{p}
p=-\frac{a\xi}{\lambda}
\end{equation}
\end{subequations}
Substituting the equalities \eqref{pq} into the first equality of the system \eqref{syst}, we obtain
that the equality \eqref{eqel} holds for \(\zeta=\xi+i\eta\). Thus we proved that
\begin{equation}
\label{incl}
(N_{\lambda}\setminus\mathbb{R})\subseteq\Gamma_{\lambda}.
\end{equation}
Let
 \begin{equation}
 \label{hp}
 \mathbb{H}^+=\lbrace\zeta:\,\im\zeta>0\rbrace,\quad \mathbb{H}^-=\lbrace\zeta:\,\im\zeta<0\rbrace
 \end{equation}
 be the upper and the lower half-plane respectively.

 According to Lemma \ref{laze}, there exist
points \(\zeta\in\mathbb{H}^{+}\setminus S\) where \(v(\zeta)>0\) and
points \(\zeta\in\mathbb{H}^{+}\setminus S\) where \(v(\zeta)<0\). This means that the set
\(N_{\lambda}\), \eqref{zers}, separates the domain \(\mathbb{H}^{+}\setminus S\). In other words,
the open set \((\mathbb{H}^{+}\setminus S)\setminus N_{\lambda}\) is disconnected.
Since \(v(\overline{\zeta})=-v(\zeta)\), the set \(N_{\lambda}\) is symmetric with respect to the real axis. The set \(\Gamma_{\lambda}\setminus N_{\lambda}\) also is symmetric with respect to the real axis.
 Since \eqref{incl}, the set \(N_{\lambda}\) can not separate the domain \((\mathbb{H}^{+}\setminus S)\)
 if \(\Gamma_{\lambda}\setminus N_{\lambda}\neq\emptyset\).\\[2.0ex]
\textsf{2}. In view of \eqref{negl}, the inequality \(0<A<B\) hold. So \(A\) is the minor semiaxis of the ellips
\(\Gamma_{\lambda}\) and \(B\) is its major semiaxis. Moreover, the inequality \(\sqrt{\frac{b}{a}}<B\)
holds. This means that the slit \(S\) is contained inside the ellipse \(\Gamma_\lambda\).
\hspace*\fill\(\Box\)\\
\begin{lemma} {\ }\\[-3.0ex]
\label{prhf}
\begin{enumerate}
\item[\textup{1.}] The functions \(u(\zeta)\) and \(v(\zeta)\) are conjugate harmonic function of \(\zeta\) in the
domain \(\zeta\in\mathbb{C}\setminus S\).
\item[\textup{2.}] The only critical points of the the functions \(u\) and \(v\) in the domain
\(\zeta\in\mathbb{C}\setminus S\) are the points
\begin{equation}
\label{Crp}
\zeta_{+}(\lambda)=\sqrt{\frac{b}{a}}\,\cdot\,\frac{|\lambda|}{\sqrt{a}}\cdot\,
\bigg(1-\frac{\lambda^2}{a}\bigg)^{-\frac{1}{2}} \ \textup{ and } \
\zeta_{-}(\lambda)=-\sqrt{\frac{b}{a}}\,\cdot\,\frac{|\lambda|}{\sqrt{a}}\cdot\,
\bigg(1-\frac{\lambda^2}{a}\bigg)^{-\frac{1}{2}},
\end{equation}
that is the points where the ellipse \(\Gamma_{\lambda}\) and the real axis \(\mathbb{R}\) intersect.
\item[\textup{3.}] If  \(\zeta\in\mathbb{H}^{+}\) lies outside the contour \(\Gamma_{\lambda}\),
then \(v(\zeta)>0\).
 If  \(\zeta\in\mathbb{H}^{+}\setminus S\) lies inside the contour \(\Gamma_{\lambda}\),
then \(v(\zeta)<0\).
\end{enumerate}
\end{lemma}
\noindent
\textsf{\textsl{Proof.}}  The functions \(u\) and \(v\) are the real and the imaginary parts of the holomorphic function
\(\sqrt{a\zeta^2+b}+\lambda\zeta\). From the Cauchy-Riemann equation it follows that the functions \(u\) and \(v\) have
the same critical points. Moreover the point \(\zeta\) is critical for \(v\) if and only if \(\zeta\) is a root of the
derivative \(a\zeta(a\zeta^2+b)^{-\frac{1}{2}}+\lambda\) of the function \(\sqrt{a\zeta^2+b}+\lambda\zeta\).
An explicit calculation shows that this derivative has only two roots \(\zeta_{+}(\lambda)\) and
\(\zeta_{-}(\lambda)\), \eqref{Crp}.

Let \(\ext(\Gamma_\lambda)\) and \(\inte(\Gamma_{\lambda})\) be the exterior and the exterior of the contour \(\Gamma_\lambda\) respectively.
Each of the sets \(\mathscr{E}^{Ext}_{\lambda}\) and \(\mathscr{E}^{Int}_{\lambda}\),
 \begin{equation}
\label{extint}
\mathscr{E}^{Ext}_{\lambda}=\ext(\Gamma_\lambda)\cap\mathbb{H}^{+},\quad
\mathscr{E}^{Int}_{\lambda}=\inte(\Gamma_\lambda)\cap(\mathbb{H}^{+}\setminus S)
\end{equation}
 is a connected open set. According to \eqref{zers} and \eqref{szs}, the continuous real valued function \(v\) does not vanish on any of these two sets. Hence the values \(v(\zeta)\) have the same sign, say \(s^{Ext}\), at all points \(\zeta\) of the set \(\mathscr{E}^{Ext}_{\lambda}\), and the same same sign, say \(s^{Int}\), at all points \(\zeta\) of the set \(\mathscr{E}^{Int}_{\lambda}\). Now the assertion 3 of Lemma \ref{prhf} is a consequence of Lemma \ref{laze}.\hspace*\fill\(\Box\)\\

\noindent
\textsf{Completion of the proof of Theorem \ref{MT}}. Let us chose the ellipse \(\Gamma_{\lambda}\)
as the contour of integration \(\Gamma\) in the integral in the right hand side of \eqref{cir}.
Since the imaginary part \(v(\zeta)\) of the exponent of the integrand
vanishes on \(\Gamma_{\lambda}\), the integral representation \eqref{cir} takes the form
\begin{equation}
\label{cirl}
\widehat{d}(\lambda,a,b)=-\frac{1}{4\pi i}\ointctrclockwise\limits_{\Gamma_\lambda}
e^{-u(\zeta)}\,d\zeta, \ \
-\sqrt{a}<\lambda<0.
\end{equation}
Since \(d\zeta=dx+idy\), we can split the integral in \eqref{cirl}:
\begin{equation}
\label{cirls}
\widehat{d}(\lambda,a,b)=-\frac{1}{4\pi i}\ointctrclockwise\limits_{\Gamma_\lambda}
e^{-u(\zeta)}\,dx(\zeta)-\frac{1}{4\pi}\ointctrclockwise\limits_{\Gamma_\lambda}
e^{-u(\zeta)}\,dy(\zeta).
\end{equation}
Since the values \(\widehat{d}(\lambda,a,b)\), \(x(\zeta)\), \(y(\zeta)\), and \(e^{-u(\zeta)}\) are real, the first integral in the right hand side of \eqref{cirls} vanishes. So the equality \eqref{cirls}
takes the form
\begin{equation}
\label{cirlt}
\widehat{d}(\lambda,a,b)=-\frac{1}{4\pi}\ointctrclockwise\limits_{\Gamma_\lambda}
e^{-u(\zeta)}\,dy(\zeta).
\end{equation}
Since the contour \(\Gamma_{\lambda}\) is symmetric with respect to the real axis \(\mathbb{R}\)
and the function \(u\) also is symmetric: \(u(\zeta)=u(\overline{\zeta})\), the equality \eqref{cirlt}
can be reduced to the form
\begin{equation}
\label{cirlr}
\widehat{d}(\lambda,a,b)=-\frac{1}{2\pi}
\sideset{}{\hspace*{-2.5ex}{\text{\raisebox{0.5ex}{\(\scriptstyle\curvearrowleft\)}}}}%
\int\limits_{\Gamma_\lambda^{+}}
e^{-u(\zeta)}\,dy(\zeta),
\end{equation}
where \(\Gamma_\lambda^{+}=\Gamma_\lambda\cap\mathbb{H}^{+}\) is the upper half of the contour
\(\Gamma_{\lambda}\). Integrating by parts in \eqref{cirlr}, we obtain
\begin{equation}
\label{cirlp}
\widehat{d}(\lambda,a,b)=\frac{1}{2\pi}
\sideset{}{\hspace*{-2.3ex}{\text{\raisebox{0.5ex}{\(\scriptstyle\curvearrowright\)}}}}%
\int\limits_{\Gamma_\lambda^{+}}
e^{-u(\zeta)}\,y(\zeta)\,du(\zeta),
\end{equation}
(The values \(y(\zeta_{\pm}(\lambda))\) at the end points \(\zeta_{\pm}(\lambda)\), \eqref{Crp},
of the integration path \(\Gamma_{\lambda}^{+}\) vanish.)

The differential \(du(\zeta)\) in \eqref{cirlp} can be represented as
\begin{equation}
\label{redi}
du(\zeta)=\frac{du(\zeta(s))}{ds}\,ds,
\end{equation}
where \(s\) is a natural parameter\footnote{\,Length of arc.
}
 on \(\Gamma_{\lambda}^{+}\).
 In other words, the differential \(du(\zeta)\) can be represented as
 \begin{equation}
\label{redv}
du(\zeta)=\frac{du}{d\vec{\tau}}(\zeta)\,ds(\zeta),
\end{equation}
 where \(\vec{\tau}(\zeta)\) is the tangent vector to the curve \(\Gamma_{\lambda}^{+}\) at the
 point \(\zeta\in\Gamma_{\lambda}^{+}\). The direction of the vector \(\vec{\tau}(\zeta)\)
 corresponds to the motion of the point \(\zeta(s)\) along the path \(\Gamma_{\lambda}^{+}\)
 from its left end point \(\zeta_{-}(\lambda)\) to the right end point \(\zeta_{+}(\lambda)\).
 If \(\vec{n}(\zeta)\) is the vector of the \emph{exterior} normal to \(\Gamma_{\lambda}^{+}\)
 at the point \(\zeta\in\Gamma_{\lambda}^{+}\), then the orientation of the frame \((\vec{\tau}(\zeta),\vec{n}(\zeta))\) coincides with the orientation of the natural frame of
 \(\mathbb{R}^2\). According the Cauchy-Riemann equations,
 \begin{equation}
 \label{CR}
 \frac{du}{d\vec{\tau}}(\zeta)=\frac{dv}{d\vec{n}}(\zeta), \ \ \forall\,\zeta\in\Gamma_{\lambda}^{+}.
 \end{equation}
 Thus the representation \eqref{cirlp} can be reduced to the form
 \begin{equation}
\label{cirlpr}
\widehat{d}(\lambda,a,b)=\frac{1}{2\pi}
\int\limits_{\Gamma_\lambda^{+}}
e^{-u(\zeta)}\,y(\zeta)\,\frac{dv}{d\vec{n}}(\zeta)\,ds(\zeta).
\end{equation}
According the assertion 3 of Lemma \ref{prhf},
\begin{equation}
\label{crpo}
\frac{dv}{d\vec{n}}(\zeta)>0, \ \ \ \forall\,\zeta\in\Gamma_{\lambda}^{+}.
\end{equation}
The inequality in \eqref{crpo} is strict because \(\frac{dv}{d\vec{n}}(\zeta)=
|\grad v(\zeta)|\) and the gradient \(\grad v(\zeta)\) of the function \(v\) vanishes
only at the critical points \(\zeta_{\pm}(\lambda)\) of the function \(v\), which are
the end points of the integration path \(\Gamma_{\lambda}^{+}\).
Evidently \(y(\zeta)>0\) and \(e^{-u(\zeta)}>0\) at every point \(\zeta\in\Gamma_{\lambda}^{+}\).
Thus the integrand
in \eqref{cirlpr} is strictly positive at every point \(\zeta\in\Gamma_{\lambda}^{+}\).
So the inequality \(\widehat{d}(\lambda,a,b)>0\) holds.
\hspace*\fill\(\Box\)

\begin{remark}
\label{prior2}
The method which we use in the second proof of \textup{Theorem \ref{MT}} is the lite version of the method which Herbert Stahl, \textup{\cite{St}}, used in his proof of the BMV conjecture.
\end{remark}

\section{The third proof of Theorem \ref{MT}.}
  For each fixed \(\eta\), the function \(\cosh\big(\eta\sqrt{t^2+\xi}\big)\) is
an entire function of the variables \(t,\xi\). Therefore, the Taylor expansion holds
\begin{equation}
\label{ser}
        \cosh\Big(\eta\sqrt{t^2+\xi}\Big)=\sum\limits_{0\leq k<\infty}\frac{1}{k!}\, \psi_{k}(t,\eta)\,\xi^{k},
\end{equation}
where
 \begin{equation*}
  \psi_{k}(t,\eta)= \frac{d^k\cosh\Big(\eta\sqrt{t^2+\xi}\Big)}{d\xi^k}\raisebox{-1.0ex}{\(\big|_{\xi=0}\)},
   \   \ \ k=0,1,2,\,\ldots\,\,.
 \end{equation*}
  It turns out that \emph{for every fixed real \(\eta\) and for every \(k=0,1,2,\,\ldots\), the function \(\psi_k(t,\eta)\) of the variable \(t\) is exponentially convex.} We prove this by induction in \(k\). Therefore for \(\xi\geq0\), the sum of the series in \eqref{ser}
  is an exponentially convex function of \(t\). To obtain Theorem \ref{MT}, we put \(\eta=\sqrt{a}\),
\(\xi=b/a\) in \eqref{ser}. (For \(a=0\), the statement of Theorem \ref{MT} is trivially true.)

Our proof of the exponential convexity of the functions \(\psi_{k}(t,\eta)\) is based on the identity
\begin{equation*}
\frac{\sinh\zeta}{\zeta}= \prod\limits_{1\leq m<\infty} \!\!\cosh\frac{\zeta}{2^m},
\end{equation*}
which holds for every \(\zeta\in\mathbb{C}\).   Substituting the expression
\begin{equation*}
\zeta= \eta\, \sqrt{t^2+\xi}
\end{equation*}
 into this identity, we obtain   the equality
\begin{equation*}
\frac{\sinh\big(\eta\sqrt{t^2+\xi}\big)}{\sqrt{t^2+\xi}} =\eta
\prod_{1\leq m<\infty} \cosh\Big(\frac{\eta}{2^m}\sqrt{t^2+\xi}\Big).
\end{equation*}
  Using the equality
 \begin{equation*}
 \dfrac{d\cosh\big(\eta\sqrt{t^2+\xi}\big)}{d\xi}=\dfrac{\eta}{2}\,
 \dfrac{\sinh\big(\eta\sqrt{t^2+\xi}\big)}{\sqrt{t^2+\xi}},
 \end{equation*}
  which holds for every \(t,\xi,\eta\),
 we obtain the equality
 \begin{equation}
 \label{Rec}
  \dfrac{d\cosh\big(\eta\sqrt{t^2+\xi}\big)}{d\xi}=
  \frac{\eta^2}{2}\prod_{1\leq m<\infty} \cosh\Big(\frac{\eta}{2^m}\sqrt{t^2+\xi}\Big).
 \end{equation}

By successive differentiation the equality \eqref{Rec} with respect to \(\xi\), we obtain
the equality
\begin{equation}
\label{sud}
\dfrac{d^{k+1}\cosh\big(\eta\sqrt{t^2+\xi}\big)}{d\xi^{k+1}}=
  \frac{\eta^2}{2}\sum\limits_{|\mathfrak{l}|=k}\bigg(\prod_{1\leq m<\infty}
  \frac{d^{\,l_{m}}\cosh\big(\frac{\eta}{2^m}\sqrt{t^2+\xi}\big)}{d\xi^{l_{m}}}\bigg),
\end{equation}
where \(k=0,1,2,3,\,\ldots\,\,\).
In \eqref{sud}, the summation is extended over all sequences\footnote{\(\textup{For} \ l_m=0, \ \ \ \dfrac{d^{\,l_{m}}\cosh\big(\frac{\eta}{2^m}\sqrt{t^2+\xi}\big)}{d\xi^{l_{m}}}
\stackrel{\textup{\tiny def}}{=}\cosh\Big(\frac{\eta}{2^m}\sqrt{t^2+\xi}\Big).\)} \(\mathfrak{l}=(l_1,l_2,l_3,\,\ldots\,)\) of non-negative integers for which \(|\mathfrak{l}|=l_1+l_2+l_3+\,\ldots\,=k.\)

The equality  \eqref{sud} holds for every \(t,\xi,\eta\).
Restricting this equality to the value \(\xi=0\), we obtain the equality
\begin{equation}
\label{sude}
\psi_{k+1}(t,\eta)=\frac{\eta^2}{2}\sum\limits_{|\mathfrak{l}|=k}\bigg(\prod_{1\leq m<\infty}
\psi_{l_{m}}\Big(t,\frac{\eta}{2^m}\Big)\bigg),
\end{equation}
which holds for every \(t\), \(\eta\), and \(k=0,1,2,3,\,\ldots\,\,\).
In \eqref{sude}, the summation is extended over all sequences
\(\mathfrak{l}=(l_1,l_2,l_3,\,\ldots\,)\) of non-negative integers for which \(|\mathfrak{l}|=l_1+l_2+l_3+\,\ldots\,=k.\)

Let \(\eta\) be an arbitrary real number. By Lemma \ref{echc}, the function
 \begin{equation}
 \label{ne}
 \psi_0(t,\eta)=\cosh \eta t
 \end{equation}
of \(t\) is exponentially convex. Moreover, the function \(\psi_0\Big(t,\dfrac{\eta}{2^m}\Big)\)
 is exponentially convex for every \(m=1,2,3,\,\ldots\,.\) (The number \(\frac{\eta}{2^m}\) here plays the
 same role as the number \(\eta\) in \eqref{ne}: it is an arbitrary real number.)

 Given \(k\geq0\), assume that all functions \(\psi_{l}(t,\frac{\eta}{2^m})\) with
 \(0\leq l\leq k\) are exponentially convex functions of \(t\).
 Then for each sequence \(\mathfrak{l}=(l_1,l_2,l_3,\,\ldots\,)\) with \(|\mathfrak{l}|=k\),
 the inequalities \(0\leq l_m\leq k\) hold. Thus, all the factors \(\psi_{l_{m}}\Big(t,\frac{\eta}{2^m}\Big)\) which appears in the product
 \(\prod\limits_{1\leq m<\infty}
\psi_{l_{m}}\Big(t,\frac{\eta}{2^m}\Big)\) are exponentially convex functions of \(t\).
Hence the product itself is an exponentially convex function. Finally, the function
\(\psi_{k+1}(t,\eta\), \eqref{sude}, which is essentially equal to the sum of all such
products with \(|\mathfrak{l}|=k\), is exponentially convex. This finishes the proof.\hspace\fill\(\Box\)
\begin{remark}
\label{coex}
Comparing the expansions  \eqref{tewb} and \eqref{ser}, we see that
\begin{equation}
\label{coeq}
\varphi_{k}(t,a)=\psi_k(t,\sqrt{a})a^{-k},\quad k=0,1,2,\,\ldots\,,\ \ t\in\mathbb{R}.
\end{equation}
As a byproduct of the third proof of \textup{Theorem \ref{MT}}, we proved that each of the functions
\(\varphi_{k}(t,a)\) is exponentially convex. Thus we have given a second proof of \textup{Theorem \ref{ptc}}.
\end{remark}
\begin{remark}
\label{Act} Actually we proved more then we formulated in Theorem \ref{MT}. Namely we proved that
for any sequence \(a_{k}(\eta)\) of non-negative numbers the sum of series
\begin{equation}
\label{sos}
s(t)=\sum\limits_{0\leq k<\infty}a_k(\eta)\psi_{k}(t,\eta)
\end{equation}
is an exponentially convex function if this series converges for every real \(t\).

If \(m\) is a positive integer and \(\xi\geq0\), then the Taylor expansion
\begin{equation}
\label{Tem}
\frac{d^{m}\cosh\big(\eta\sqrt{t^2+\xi})}{d\xi^{m}}=
\sum\limits_{m\leq k<\infty}\frac{1}{(k-m)!}\varphi_{k}(t,\eta)\,\xi^{k-m}
\end{equation}
is of the form \eqref{sos} with \(a_{k}(\eta)=0\) for \(0\leq k <m\), \(a_{k}(\eta)=\dfrac{1}{(k-m)!}\xi^{k-m}\) for \(k=m,m+1,m+2,\,\ldots\,\,.\)
\end{remark}

In particular, for \(m=1\) the following result holds:
\begin{theorem}
\label{meo}
For any \(a>0\) and \(b>0\), the function \(\psi(t)=\dfrac{\sinh\big(\sqrt{at^2+b}\big)}{\sqrt{at^2+b}}\)
is an exponentially convex function of the variable \(t\).
\end{theorem}

\end{document}